\documentclass[12pt]{amsart}
\usepackage{amssymb, amsmath, amscd, amsthm}
\usepackage[margin=1.1in]{geometry}
\usepackage{lineno}
\usepackage{amsfonts}
\usepackage{psfrag}
\usepackage[all]{xy}
\usepackage{graphicx}
\usepackage{verbatim}
\usepackage{color}
\usepackage{tikz-cd}
\usepackage{enumerate}
\usepackage[shortlabels]{enumitem}

\newtheorem{thm}{Theorem}[section]

\newtheorem{prop}[thm]{Proposition}
\theoremstyle{definition}
\newtheorem{defn}[thm]{Definition}
\theoremstyle{remark}

\numberwithin{equation}{subsection}

\theoremstyle{definition}

\newtheorem*{thm*}{Theorem}

\newcommand{\dgm}{\mathrm{Dgm}}

\title[Geometry of the matching distance]{Geometry of the matching distance for $2$D filtering functions
}
\author[Ethier]{Marc Ethier}
\email{Marc.Ethier@uqat.ca}
\address{}
\author[Frosini]{Patrizio Frosini }
\email{patrizio.frosini@unibo.it}
\address{}
\author[Quercioli]{Nicola Quercioli}
\email{nicola.quercioli@enea.it}
\address{}
\author[Tombari]{Francesca Tombari}
\email{tombari@kth.se}
\address{}


\begin{document}
\maketitle

\begin{abstract}
In this paper we exploit the concept of extended Pareto grid to study the geometric properties of the matching distance for $\mathbb{R}^2$-valued  regular functions defined on a Riemannian closed manifold. In particular, we prove that in this case the matching distance is realised either at special values or at values corresponding to vertical, horizontal or slope $1$ lines.
\end{abstract}

\section{Introduction}\label{sec1}

Feature extraction and comparison are the main tasks of data analysis.
In topological data analysis this translates into the problem of comparing persistence modules, which encode the homological features extracted from geometric objects. 
In order to be able to compare persistence modules a distance is needed. 
There is a wide variety of distances in the space of $1$-parameter persistence modules, 
such as the bottleneck and Wasserstein distances.
However, such distances do not directly generalise to the multiparameter setting.
Thus, different ones have been proposed over the past years, turning into a substantial catalogue, 
see for example~\cite{pers_Betti_stable, coher_match, inter}. 
One of those is the matching distance, which can be defined in particular for 2-parameter persistence modules. 
This pseudometric, introduced in~\cite{pers_Betti_stable}, 
is a generalisation of the classical bottleneck distance for $1$-parameter persistence modules and measures the difference between the $2$-dimensional Betti numbers functions (also known as rank invariants) of persistence modules.
The definition of matching distance is based on a foliation method, consisting of ``slicing'' the $2$-parameter persistence module into infinitely many $1$-dimensional components
by means of lines of positive slope, that we refer to as filtering lines.
The matching distance is then obtained by taking the supremum over all such lines of the bottleneck distances between the resulting persistence diagrams after a suitable normalisation. 

According to the definition, in order to compute the matching distance between persistence modules,
one should take into account infinitely many bottleneck distance computations. 
Many efforts have been devoted to make this computation efficient, for example, by identifying a finite number of filtering lines contributing to the actual computation~\cite{comp_disc, comp_asymp, comp_exact} or approximation~\cite{comp1,comp2, efficient-approx} of the distance.
However, what many of these works have in common is that their starting point is a pair of $2$-parameter persistence modules. 
Our approach is similar in scope, but different in nature.
We consider regular filtering functions on a smooth manifold with values in $\mathbb{R}^2$. 
Their sublevel-set filtrations still return $2$-parameter persistence modules, 
for which we can compute the $2$-dimensional persistent Betti numbers functions and, hence, the matching distance between them. 
Our approach allows us to observe phenomena and exploit structures that are not visible when directly considering persistence modules.  
For example, it is possible to exploit the differentiable structure of the filtering functions to
identify points in the persistence diagrams associated to each filtering line. 
This structure made of arcs and half-lines is known as extended Pareto grid~\cite{coher_match} (see also~\cite{wan}).
The convenience of such an approach relies on the fact that, by using this approach, the changes in homology that occur when the filtering line changes are easy to follow and control. 

In this context of $2$-parameter persistence modules derived from regular filtering functions on smooth manifolds, 
we show that filtering lines of slope $1$ play a special role in the computation of the matching distance.
Our main result shows that the matching distance between the $2$-dimensional persistent Betti numbers functions of two filtering functions is indeed
realised either on values corresponding to vertical, horizontal or slope $1$ lines, or on special values associated with the two functions. 
The authors of~\cite{comp_disc} recently obtained an analogous result in the discrete setting. They show that the matching distance is realised either on values corresponding to diagonal lines or on what they call switch values. 
One main difference is that the collection of special values that we encounter, called special set, is strongly related with the differentiable structure of our input. 
In particular, it relies on the structure of the extended Pareto grid associated with a function and on the Position Theorem, proved in~\cite{coher_match},
which relates points of a persistence diagram to points in the extended Pareto grid. 

In this paper, we aim to prove the following:
\begin{thm*}
The matching distance between $f$ and $g$ is realised either on a value associated with a line of slope $1$,  a vertical or horizontal line, or on a special value of $(f,g)$.
\end{thm*}

\section{Matching distance}\label{sec:setting}
Let $M$ be a closed $C^{\infty}$-manifold with a Riemannian metric defined on it. 
Let $f\colon M\to \mathbb{R}$ be a smooth function. The filtered homology of the sublevel sets of $f$ is known as persistent homology. 
This information can be encoded as a multiset of points in $\{(x,y)\in\mathbb{R}^2\mid x\le y\}$, known as the persistence diagram of $f$ and denoted by $\dgm(f)$. 
The subset $\Delta =\{(x,y)\in \mathbb{R}^2\mid x= y\}$ is always considered to be in the persistence diagram of a function and, 
by convention, we treat it as a unique point with infinite multiplicity. 
See~\cite{filtered-persistence} for more details about $1$-parameter persistent homology for sublevel set filtrations.

\begin{figure}
    \centering
    \includegraphics[width=10cm]{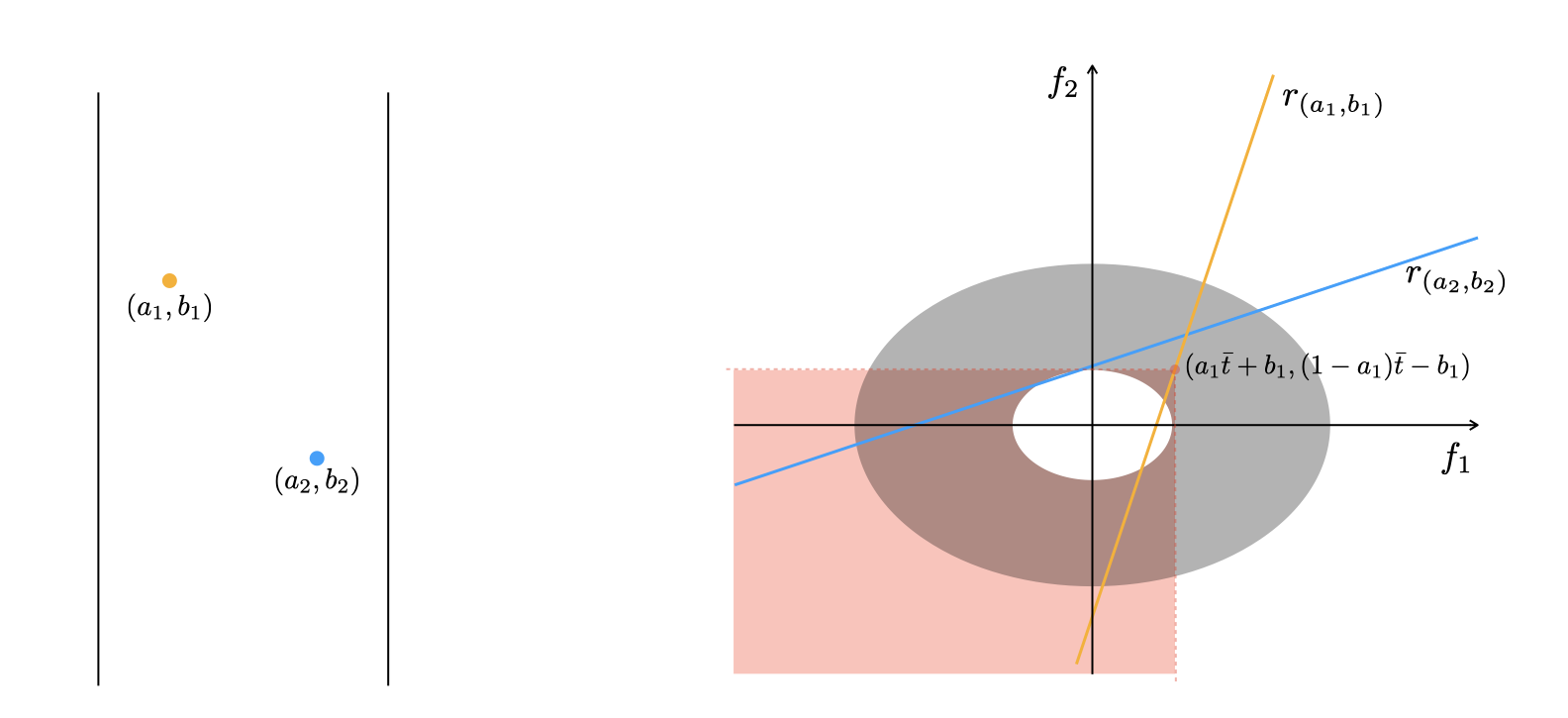}
    \caption{On the left the space of parameters $]0,1[\times\mathbb{R}$ is represented. On the right we show the filtering lines corresponding to the two parameters selected and the projection of a torus on the plane. Each of the lines gives a sublevel set filtration of the torus in the direction of the line. 
    The area down and left of the red point $(a_1\bar t+b_1, (1-a_1)\bar t-b_1)$ on the line $r_{(a_1,b_1)}$ is the sublevel set $M_{\bar t}^{a_1,b_1}$ associated to such a point. }
    \label{foliation}
\end{figure}

Let $f=(f_1,f_2)\colon M\to \mathbb{R}^2$ and $g=(g_1,g_2)\colon M\to \mathbb{R}^2$ be smooth functions. Consider the set of pairs $(a,b)$
in $]0,1[\times \mathbb{R}$ with the uniform metric $d_{\infty}((a,b), (a',b'))=\max\{\lvert a-a'\rvert, \lvert b-b'\rvert\}$. It parameterises all the lines of $\mathbb{R}^2$ with positive slope in the following way: $r_{(a,b)}$ is the line containing points of the form $t(a,1-a)+(b,-b)$ with $t$ in $\mathbb{R}$. Each point $(u(t),v(t))$ of $r_{(a,b)}$ can be associated with the set $M_t^{a,b}=M_{(u(t),v(t))}=\allowbreak \left\{x\in M \mid f_1(x)\le u(t) \text{ and } f_2(x)\le v(t)\right\}$. 
This defines a $1$-dimensional filtration, depending on the line $r_{(a,b)}$, which can be associated with a persistence diagram. 
Letting $(a,b)$ vary, one obtains a collection of persistence diagrams described by the $2$D persistent Betti numbers function of $f$. 
As observed in~\cite{pers_Betti_stable}, $M_{t}^{(a,b)}$ is also equal to $\left\{x\in M\mid f_{(a,b)}(x)\le t\right\}$, where $f_{(a,b)}(x)=\max\left\{\frac{f_1(x)-b}{a}, \frac{f_2(x)+b}{1-a}\right\}$ (see Figure~\ref{foliation}). However, more commonly, $f_{(a,b)}$ is normalised, without changing the nature of the filtration, and $f^*_{(a,b)}=\min\{a,1-a\}f_{(a,b)}$ is instead considered.    

Given two functions $f$ and $g$, the  \emph{matching distance}~\cite{pers_Betti_stable} is defined by 
\[
D_{\text{match}}(f,g)=\sup _{(a,b)\in ]0,1[\times \mathbb{R}}d_B\left(\dgm \left(f^*_{(a,b)}\right), \dgm \left(g^*_{(a,b)}\right)\right).
\]
Here, $d_B\left(\dgm \left(f^*_{(a,b)}\right), \dgm \left(g^*_{(a,b)}\right)\right)$ is the bottleneck distance between the persistence diagrams of $f^*_{(a,b)}$ and $g^*_{(a,b)}$, i.e.,
\[
d_B\left(\dgm \left(f^*_{(a,b)}\right), \dgm \left(g^*_{(a,b)}\right)\right)=\inf _{\sigma}\text{cost}(\sigma),
\]
where $\text{cost}(\sigma)=\max_{X\in \dgm\left(f^*_{(a,b)}\right)} d\left(X,\sigma (X)\right)$, $\sigma$ runs over all bijections, called matchings, between $\dgm\left(f^*_{(a,b)}\right)$ and $\dgm\left(g^*_{(a,b)}\right)$ and, if $X=(x_1,x_2)$ and $Y=(y_1,y_2)$, then
\[
d(X,Y) =
\begin{cases}
\kappa & \textup{if } X=(x_1,x_2),\; Y=(y_1,y_2) \in \Delta^+, \\
\vert x_1 - y_1 \vert & \textup{if } X = (x_1,\infty),\; Y = (y_1,\infty), \\
\frac{x_2-x_1}{2} & \textup{if } X = (x_1,x_2) \in \Delta^+, \; Y = \Delta, \\
\frac{y_2-y_1}{2} & \textup{if } Y = (y_1,y_2) \in \Delta^+, \; X = \Delta, \\
0 & \textup{if } X = Y = \Delta, \\
\infty & \textup{otherwise.}
\end{cases}
\]
where $\kappa=\min \{ \max \{ \vert x_1 - y_1 \vert, \vert x_2 - y_2 \vert\},\max \{ \frac{x_2-x_1}{2}, \frac{y_2-y_1}{2} \}\}$.
A matching realising the matching distance, whenever it exists, is called an \emph{optimal matching}.
Note that the matching distance $D_{\text{match}}(f,g)$ can be seen both as a pseudo-metric between the persistent Betti numbers functions of $f$ and $g$, and between the filtering functions $f$ and $g$. For the sake of simplicity, we keep the notation $D_{\text{match}}(f,g)$ for both cases. For more details about this definition and the foliation method we refer to~\cite{pers_Betti_stable}. 

\section{Extended Pareto grid}
In this section we recall the relation between a differential construction associated with a smooth function $f\colon M\to \mathbb{R}^2$, 
called the extended Pareto grid, and the points of the persistence diagrams $\dgm\left(f^*_{(a,b)}\right)$.
This connection is established in the Position Theorem proved in~\cite{coher_match}. 

Recall that the \emph{Jacobi set} of $f$ is the collection 
\[
\mathbb{J}(f)=\{p\in M\mid \nabla f_1=\lambda \nabla f_2 \text{ or }\nabla f_2=\lambda \nabla f_1, \text{ for some }\lambda\in \mathbb{R}\}.
\]
The \emph{Pareto critical set} of $f$ is the subset of $\mathbb{J}(f)$ given by
\[
\mathbb{J}_P(f)=\left\{p\in \mathbb{J}(f)\mid \nabla f_1=\lambda \nabla f_2 \text{ or }\nabla f_2=\lambda \nabla f_1, \text{ for some }\lambda\le 0\right\}.
\]

Assume now that $f$ is not only smooth, but it also satisfies the following properties:
\begin{enumerate}[(i)]
\item No point $p$ exists in $M$ at which both $\nabla f_1$ and $\nabla f_2$ vanish.
\item $\mathbb{J}(f)$ is a $1$-manifold smoothly embedded in $M$ consisting of finitely many components, each one diffeomorphic to a circle.
\item $\mathbb{J}_P(f)$ is a $1$-dimensional closed submanifold of M, with boundary in $\mathbb{J}(f)$.
\item If we denote by $\mathbb{J}_C(f)$ the subset of $\mathbb{J}(f)$ where $\nabla f_1$ and $\nabla f_2$ are orthogonal to $\mathbb{J}(f)$, then the connected components of $\mathbb{J}_P(f)\setminus\mathbb{J}_C(f)$ are finite in number, each one being
diffeomorphic to an interval. With respect to any parameterisation of each component, one of $f_1$ and $f_2$ is strictly increasing and the other is strictly decreasing.
Each component can meet critical points for $f_1$, $f_2$ only at its endpoints.
\end{enumerate}

\begin{figure}
    \centering    \includegraphics[width=10cm]{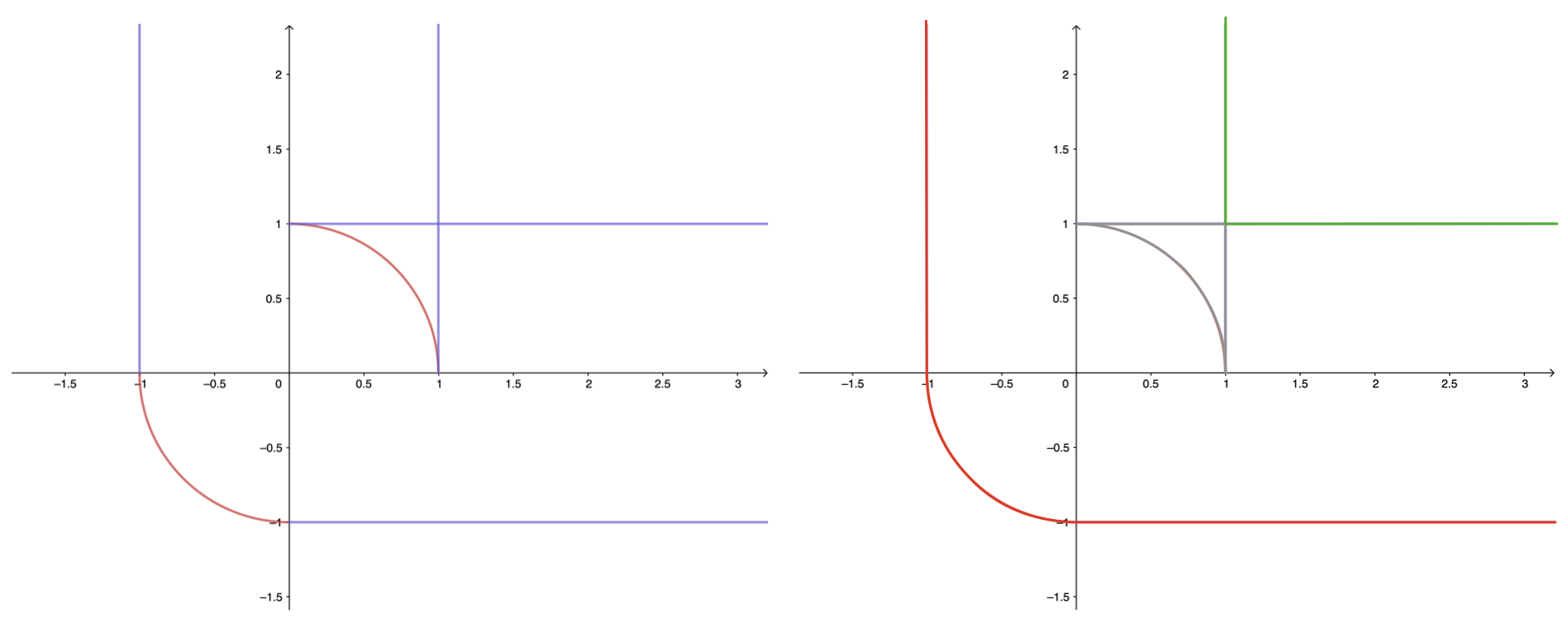}
    \caption{The extended Pareto grid of the function $f(x,y,z)=(x,z)$ on $S^2=\{(x,y,z)\in \mathbb{R}^3\mid x^2+y^2+z^2=1\}$.On the left, the red arcs are the proper contours and the violet half-lines are the improper contours. 
    Any point in this extended Pareto grid corresponds to a birth or death of a homological class. For example, on the right side of this figure, the red points correspond to the birth of homology classes in degree $0$, while the green points correspond to the birth of homology classes in degree $2$.}
    \label{ghsjdkucdde}
\end{figure}

Denote by $\{p_1,\dots, p_h\}$ and $\{q_1,\dots, q_k\}$, respectively, the critical points of $f_1$ and $f_2$. Since the function $f$ satisfies (i), then $\{p_1,\dots, p_h\}\cap\{q_1,\dots, q_k\}=\emptyset $. The \emph{extended Pareto grid} of $f$ is defined as the union 
\[
\Gamma (f)=f\left(\mathbb{J}_P(f)\right)\cup \left(\bigcup _i v_i\right) \cup  \left(\bigcup _j h_j\right)
\]
where $v_i$ is the vertical half-line $\{(x,y)\in \mathbb{R}^2\mid x=f_1(p_i), y\ge f_2(p_i)\}$ and
$h_j$ is the horizontal half-line $\{(x,y)\in \mathbb{R}^2\mid x\ge f_1(q_j), y= f_2(q_j)\}$.
We refer to these half-lines as \emph{improper contours} and to the closure of the image of the connected components of $\mathbb{J}_P(f)\setminus\mathbb{J}_C(f)$ as \emph{proper contours} of $\Gamma (f)$. 
Figure~\ref{ghsjdkucdde} shows an example of extended Pareto grid for the projection of a sphere in $\mathbb{R} ^3$ on the plane $y=0$. 
The violet horizontal half-lines originate at critical values of $f_2$, while the vertical ones originate at critical values of $f_1$. 
The red arcs are the images of those arcs on the sphere in which the gradients $\nabla f_1$ and $\nabla f_2$ have the same direction but opposite orientation.
Observe that, because of property (ii), the number of contours in $\Gamma (f)$ is finite.
Moreover, property (iv) ensures that every contour can be parameterised as a curve whose two coordinates are respectively strictly decreasing and strictly increasing.
For more details about properties (i)-(iv) we refer the interested reader to~\cite{coher_match, wan}.

One may observe that the portions of contours delimited by points of intersection between different contours correspond to births and deaths of homology classes. 
For example, the red union of contours corresponds to the birth of a homology class in degree 0 and the green portions of contour to the birth of a homology class in degree 2. 
For a richer example we refer the reader to~\cite[Figure 8]{coher_match}.

The Position Theorem (Theorem 2 in~\cite{coher_match}) allows us to obtain the coordinates of the points in the persistence diagram of $f^*_{(a,b)}$ just by looking at the extended Pareto grid of the function and the filtering line $r_{(a,b)}$.
It reads as follows:
\begin{thm}\label{posth}
Let $(a,b)$ be in $]0,1[\times\mathbb{R}$ and $X$ in $\dgm (f^*_{(a,b)})\setminus \{\Delta\}$. Then, for each finite coordinate $w$ of $X$, a point $(p_1,p_2)$ in $r_{(a,b)}\cap \Gamma(f)$ exists such that $w=\frac{\min\{a, 1-a\}}{a}(p_1-b)=\frac{\min\{a, 1-a\}}{1-a}(p_2+b)$.
\end{thm}

In~\cite{coher_match} the set of filtering functions considered is the set of normal functions. However, the reader can observe that the proof of this specific theorem is actually independent from this assumption and it is valid also in our current setting.

\section{Extension of persistence diagrams}\label{sec:extension}
In this section we show that it is possible to extend each $2$D persistent Betti numbers function from the open set $]0,1[\times \mathbb{R}$, where it is defined, to the closed set $[0,1]\times \mathbb{R}$.
Moreover, we prove that the matching distance between $f$ and $g$ can be realised on the compact set $[0,1]\times [-\overline{C}, \overline{C}]$, with $\overline{C}=\max \{\lVert f\rVert_{\infty}, \lVert g\rVert_{\infty}\}$.

\begin{prop}\label{propineq1}
Let $C$ be a positive real number. If $0<a,a'\le \frac{1}{2}$ and $\lvert b\rvert\le C$, then, for every $b'$,
\[
\left\lVert f^*_{(a,b)}-f^*_{(a',b')}\right\rVert_\infty\le 4\lvert a-a'\rvert \left(\lVert f_2\rVert_\infty+C\right)+ 3\lvert b-b'\rvert.
\]
\end{prop}

\begin{proof}
Since $a,a'\le \frac{1}{2}$, then $\min\left\{a,1-a\right\}=a$ and $\min\left\{a',1-a'\right\}=a'$. 
Therefore, recalling that $\lvert \max \left\{\alpha,\beta\right\}-\max\{\gamma,\delta\}\rvert\le \max\{\lvert\alpha-\gamma\rvert,\lvert\beta-\delta\rvert\}$ and observing that $(1-a)(1-a')\ge \frac{1}{4}$,
\begingroup
\allowdisplaybreaks
\begin{align*}
& \left\lVert f^*_{(a,b)}-f^*_{(a',b')}\right\rVert_\infty = 
\left\lVert a\max \left\{\frac{f_1-b}{a},  \frac{f_2+b}{1-a}\right\}-
a' \max \left\{\frac{f_1-b'}{a'},\frac{f_2+b'}{1-a'}\right\}
\right\rVert_\infty \nonumber\\
&= \sup_{x\in M} \left\lvert
\max \left\{f_1(x)-b,  \frac{a}{1-a}(f_2(x)+b)\right\}  - \max \left\{f_1(x)-b', \frac{a'}{1-a'}(f_2(x)+b')\right\} \right\rvert \nonumber\\
&\le \sup_{x\in M} 
\max \left\{\lvert b-b'\rvert,  \left\lvert\frac{a}{1-a}(f_2(x)+b)-\frac{a'}{1-a'}(f_2(x)+b')\right\rvert\right\} \nonumber\\
&=\sup_{x\in M}
\max \left\{\lvert b-b'\rvert,  \frac{\left\lvert af_2(x)+ab -aa'b-a'f_2(x)-a'b'+aa'b'\right\rvert}{(1-a)(1-a')}\right\}
\nonumber\\
&\le \max \left\{\lvert b-b'\rvert,  \frac{\lvert a-a'\rvert\lVert f_2\rVert_\infty+\lvert ab-a'b'\rvert + \lvert b-b'\rvert aa'}{(1-a)(1-a')}\right\} \nonumber\\
&\le
\max \left\{\lvert b-b'\rvert ,  4\lvert a-a'\rvert \lVert f_2\rVert_\infty+4\lvert ab-a'b'\rvert+ 4\lvert b-b'\rvert aa'\right\} \nonumber\\
&\le
\max \{\lvert b-b'\rvert ,  4\lvert a-a'\rvert \lVert f_2\rVert_\infty+4\lvert ab-a'b'\rvert+ \lvert b-b'\rvert\} \nonumber\\
&=
4\lvert a-a'\rvert \lVert f_2\rVert_\infty+4\lvert ab-a'b'\rvert + \lvert b-b'\rvert \nonumber\\
&\le
4\lvert a-a'\rvert\lVert f_2\rVert_\infty+4\lvert a-a'\rvert \lvert b\rvert +4\lvert b-b'\rvert a'+ \lvert b-b'\rvert \nonumber\\
&\le
4\lvert a-a'\rvert\lVert f_2\rVert_\infty+4\lvert a-a'\rvert\lvert b\rvert+ 3\lvert b-b'\rvert \nonumber\\
&\le 
4\lvert a-a'\rvert(\lVert f_2\rVert_\infty+C)+ 3\lvert b-b'\rvert. \nonumber
\end{align*}
\endgroup
\end{proof}
By observing that, if $f=(f_1,f_2)$ and $h=(f_2,f_1)$, then $f^*_{(a,b)}=h^*_{(1-a,-b)}$,
we obtain an analogous result to Proposition~\ref{propineq1} for $a,a'\ge \frac{1}{2}$.

\begin{prop}\label{sertyghvcxd}
Let $C$ be a positive real number. If $\frac{1}{2}\le a,a'< 1$ and $\lvert b\rvert\le C$, then, for every $b'$,
\[
\left\lVert f^*_{(a,b)}-f^*_{(a',b')}\right\rVert_\infty\le 4\lvert a-a'\rvert (\lVert f_1\rVert_\infty+C)+ 3\lvert b-b'\rvert.
\]
\end{prop}

As a consequence, the function
\[
f^*_{(\cdot,\cdot)}\colon ]0,1[\times \mathbb{R} \to C(M,\mathbb{R}) , \quad (a,b) \mapsto f^*_{(a,b)}
\]
is locally Lipschitz. This is the content of the following result:

\begin{thm}\label{ghjywuikdc}
If $\lvert b\rvert\le C$, then for every $0<a,a'<1$ and every $b'$,
\[
\left\lVert f^*_{(a,b)}-f^*_{(a',b')}\right\rVert_{\infty}\le 4\lvert a-a'\rvert (\lVert f\rVert_\infty+C)+ 3\lvert b-b'\rvert.
\]
\end{thm}
\begin{proof}
If $a,a'\le \frac{1}{2}$ or $a,a'\ge \frac{1}{2}$, the statement follows directly from Propositions \ref{propineq1} and \ref{sertyghvcxd}.
Without loss of generality, we can assume that $a \le \frac{1}{2}$ and $a' \ge \frac{1}{2}$. Moreover, consider $\left(\frac{1}{2},b\right), \left(\frac{1}{2},b'\right)$.
We have that
\begin{align*}
     \left\lVert f^*_{(a,b)} - f^*_{(a',b')} \right\rVert_\infty & \le  \left\lVert f^*_{(a,b)} - f^*_{(\frac{1}{2},b)} \right\rVert_\infty + \left\lVert f^*_{(\frac{1}{2},b)} - f^*_{(\frac{1}{2},b')}\right\rVert_\infty + \\ &+ \left\lVert f^*_{(\frac{1}{2},b')} - f^*_{(a',b')} \right\rVert_\infty \\ 
    & \le 4\left\lvert a-\frac{1}{2} \right\rvert \left(\lVert f\rVert_\infty+C\right)+ 3\lvert b-b\rvert + 4\left\lvert \frac{1}{2}- \frac{1}{2}\right\rvert (\lVert f\rVert_\infty+C)+ \\
    & + 3\lvert b-b'\rvert  + 4\left\lvert \frac{1}{2}-a'\right\rvert (\lVert f\rVert_\infty+C)+ 3\lvert b'-b'\rvert \\
    & = 4 \left(\left\lvert a-\frac{1}{2}\right\rvert + \left\lvert \frac{1}{2}-a'\right\rvert\right) (\lVert f\rVert_\infty+C) + 3\lvert b-b'\rvert  \\ & = 4\lvert a-a'\rvert (\lVert f\rVert_\infty+C)+ 3\lvert b-b'\rvert.
\end{align*}
\end{proof}
 
In Theorem \ref{ghjywuikdc} we showed that the function $f^*_{(\cdot,\cdot)}$ is locally Lipschitz. As such, it can be extended to the parameter values $(0,b)$ (resp. $(1,b)$) as the limit $f^*_{(0,b)}:=\lim_{(a',b')\to (0,b)}f^*_{(a',b')}$ 
(resp., $f^*_{(1,b)}:=\lim_{(a',b')\to (1,b)}f^*_{(a',b')}$), for every $b$ in $\mathbb{R}$. 
Such a function is continuous, 
and the stability of persistence diagrams with respect to the uniform norm implies that the limit
$\lim_{(a',b')\to (0,b)}\dgm\left(f^*_{(a',b')}\right)$ (resp. $\lim_{(a',b')\to (1,b)}\dgm\left(f^*_{(a',b')}\right)$) also exists and is equal to 
$\dgm\left(f^*_{(0,b)}\right)$ (resp. $\dgm\left(f^*_{(1,b)}\right)$).
In other words, $f^*_{(\cdot,\cdot)}$ can be uniquely extended to $[0,1] \times \mathbb{R}$ and this extension is also a locally Lipschitz function.
Therefore, in the rest of this paper, we will be allowed to consider the functions $f^*_{(a,b)}$ for any 
$(a,b)$ in $[0,1]\times\mathbb{R}$. 
The limit functions $f_{(0,b)}^*$ and $f^*_{(1,b)}$ can be computed explicitly for any $b$ in $\mathbb{R}$: 
\begin{align*}
f^*_{(0,b)}(x) & =\lim _{(a',b')\to (0,b)}f^*_{(a',b')}(x) \\ &=\lim _{a'\to 0} a'\max\left\{ \frac{f_1(x)-b}{a'}, \frac{f_2(x)+b}{1-a'}\right\}\\&=\max\left\{ f_1(x)-b,0\right\}, 
\end{align*}
\begin{align*}
f^*_{(1,b)}(x) & =\lim _{(a',b')\to (1,b)}f^*_{(a',b')}(x)\\&= \lim _{a'\to 1} (1-a')\max\left\{ \frac{f_1(x)-b}{a'}, \frac{f_2(x)+b}{1-a'}\right\}\\&=\max\left\{0, f_2(x)+b\right\}.
\end{align*}

Since Theorem~\ref{ghjywuikdc} enables us to extend the functions $f^*_{(\cdot, \cdot)}$ and $g^*_{(\cdot,\cdot)}$ to $[0,1]\times \mathbb{R}$, the function
\[
(a,b)\mapsto d_B\left(\dgm\left(f^*_{(a,b)}\right), \dgm\left(g^*_{(a,b)}\right)\right) \tag{4.1} \label{contdb}
\]
can be extended to $[0,1]\times \mathbb{R}$, too. Furthermore, it is continuous because of the stability of persistence diagrams and, hence, it admits a maximum in its compact domain. 

Next, we show that it is not restrictive to compute the matching distance for parameters in $[0,1]\times [-\overline{C},\overline{C}]$, where $\overline{C}=\max \{\lVert f\rVert_{\infty},\lVert g\rVert_{\infty}\}$. 

\begin{prop}\label{wertyuiolkmnbvfrtyh}
There exists $(\bar a,\bar b)$ in $[0,1]\times [-\overline{C},\overline{C}]$, with $\overline{C}=\max \{\lVert f\rVert_{\infty},\lVert g\rVert_{\infty}\}$, such that 
\begin{align*}
D_{\mathrm{match}}(f,g)& =\max_{(a,b)\in [0,1]\times [-\overline{C}, \overline{C}]} d_B\left(\dgm\left(f_{( a,b)}^*\right),\dgm\left(g_{(a, b)}^*\right)\right)
\\& =d_B\left(\dgm\left(f_{(\bar a,\bar b)}^*\right),\dgm\left(g_{(\bar a,\bar b)}^*\right)\right).
\end{align*}

\end{prop}

\begin{proof}

Our strategy is to check what happens when $\lvert b\rvert \ge \overline{C}$. There are four possible cases given by the combinations of $a\le \frac{1}{2}$ or $a\ge \frac{1}{2}$ and $b\le -\overline{C}$ or $b\ge \overline{C}$. 
Consider the case $a\le \frac{1}{2}$ and $b\le -\overline{C}$. We have $f^*_{(a,b)}=a\max \left\{\frac{1}{a}(f_1-b), \frac{1}{1-a}(f_2+b)\right\}$. 
However, $\frac{1}{a}(f_1-b)\ge \frac{1}{a}(f_1+\overline{C})\ge 0$ and $\frac{1}{1-a}(f_2+b)\le \frac{1}{1-a}(f_2-\overline{C})\le 0$. Thus, $f^*_{(a,b)}=f_1-b$ and, similarly, $g^*_{(a,b)}=g_1-b$. 
The bottleneck distance between their persistence diagrams will thus be $d_B\left(\dgm\left(f^*_{(a,b)}\right),\dgm\left(g^*_{(a,b)}\right)\right)=\allowbreak d_B\left(\dgm\left(f_1-b\right),\dgm\left(g_1-b\right)\right)=\allowbreak d_B\left(\dgm\left(f_1\right),\dgm\left(g_1\right)\right)$. 
Therefore, $d_B\left(\dgm\left(f^*_{(a,b)}\right),\dgm\left(g^*_{(a,b)}\right)\right)$ is constant for $a\le \frac{1}{2}$ and $b\le -\overline{C}$. 
Hence we can limit ourselves to computing its value for $a= \frac{1}{2}$ and $b=-\overline{C}$.

Consider now $a\ge \frac{1}{2}$ and $b\le -\overline{C}$. We have $f^*_{(a,b)}=\frac{1-a}{a}(f_1-b)$ and, similarly, $g^*_{(a,b)}=\frac{1-a}{a}(g_1-b)$. Fixing $a$, we observe that in this case $d_B\left(\dgm\left(f^*_{(a,b)}\right),\dgm\left(g^*_{(a,b)}\right)\right)=\frac{1-a}{a}d_B\left(\dgm(f_1-b),\dgm(g_1-b)\right)=\frac{1-a}{a}d_B(\dgm(f_1),\dgm(g_1))$ is constant with respect to $b$. 
Since $a\ge\frac{1}{2}$ was chosen arbitrarily, and there is no dependence on $b$, we can choose them to be $a= \frac{1}{2}$ and $b=-\overline{C}$ and conclude. 

The other two cases follow the same strategy. 
\end{proof}
\begin{figure}
    \centering
    \includegraphics[width=7cm]{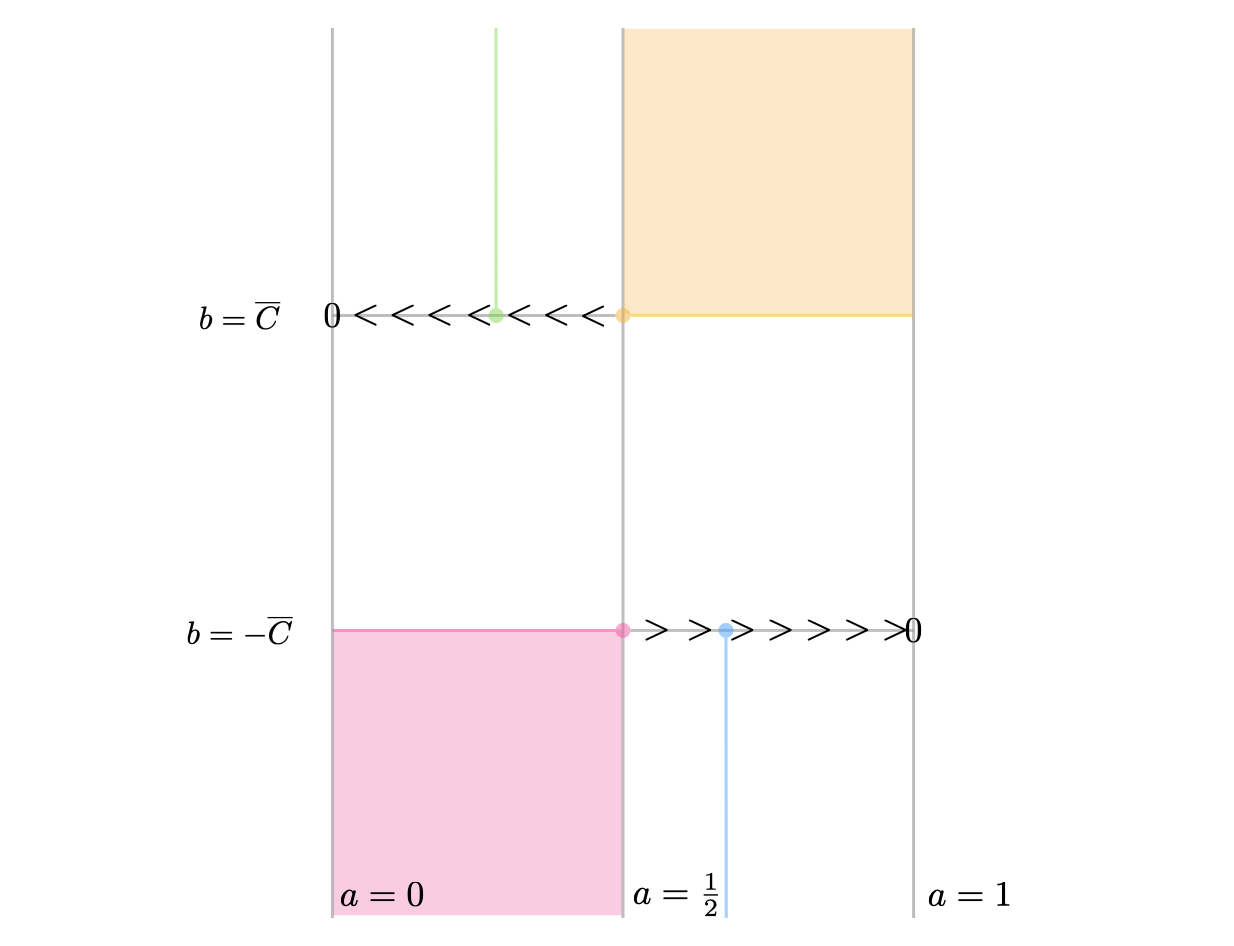}
    \caption{The bottleneck distance is constant on the coloured regions and half-lines, whereas it is non-increasing with respect to $a$ on $b=-\overline{C}$, if $a\ge \frac{1}{2}$, and non-decreasing with respect to $a$ on $b=\overline{C}$, if $a\le \frac{1}{2}$.}
    \label{box}
\end{figure}

The above proof also shows that the continuous function $d_B\left(\dgm\left(f^*_{(a,b)}\right), \dgm\left(g^*_{(a,b)}\right)\right)$ is constant on the
segments $\{(a,b)\mid 0\le a\le\frac{1}{2},
b=-\overline{C}\}$ and $\{(a,b)\mid \frac{1}{2} \le a \le 1, b=\overline{C}\}$, non-increasing on the segment $\{(a,b)\mid 1\ge a\ge\frac{1}{2},
b=-\overline{C}\}$ and non-decreasing on the segment $\{(a,b)\mid 0\le a\le\frac{1}{2}, b=\overline{C}\}$. 
Moreover, it is $0$ on $(0, \overline{C})$ and $(1,-\overline{C})$ (see Figure~\ref{box}).
Furthermore, we would like to point out that Proposition \ref{wertyuiolkmnbvfrtyh} gives us a new formulation for the definition of the matching
distance $D_{\mathrm{match}}$ as follows:
\begin{align*}
D_{\mathrm{match}}(f,g)& =\max_{(a,b)\in [0,1]\times \mathbb{R}} d_B\left(\dgm\left(f_{( a,b)}^*\right),\dgm\left(g_{(a, b)}^*\right)\right).
\end{align*}

\section{Special set and matching distance}

In this section we introduce the special set associated with a pair of functions $(f,g)$. We prove that the matching distance between two functions is realised either on values associated with vertical, horizontal or slope $1$ lines, or on this special set.

\begin{defn}
Let $\text{Ctr}(f,g)$ be the set of all curves that are contours of $f$ or $g$.
The \emph{special set}
of $(f,g)$, denoted by $\text{Sp}(f,g)$, is the
collection of all $(a,b)$ in $]0,1[\times [- \overline{C},
\overline{C}]$ for which two distinct pairs $\{\alpha_p, \alpha_q\}$,
$\{\alpha_s, \alpha_t\}$ of contours in $\text{Ctr}(f,g)$ intersecting $r_{(a,b)}$
exist, such that $\{\alpha_p, \alpha_q\} \ne \{\alpha_s, \alpha_t\}$ and
\begin{itemize}
\item $c_1\lvert x_P-x_Q\rvert=c_2\lvert x_{S}-x_{T}\rvert$, with $c_1,c_2\in \{1,2\}$, if $a\le \frac{1}{2}$,
\item $c_1\lvert y_P-y_Q\rvert=c_2\lvert y_{S}-y_{T}\rvert$, with $c_1,c_2\in \{1,2\}$, if $a\ge \frac{1}{2}$,
\end{itemize}
where $P=P_{(a,b)}=r_{(a,b)}\cap \alpha_p$,
$Q=Q_{(a,b)}=r_{(a,b)}\cap \alpha_q$, $S=S_{(a,b)}=r_{(a,b)}\cap \alpha_s$ and
$T=T_{(a,b)}=r_{(a,b)}\cap \alpha_t$, and $x_{*}$, $y_{*}$ denote abscissas
and ordinates of these points.
An element of the special set $\text{Sp}(f,g)$ is called a
\emph{special value} of the pair $(f,g)$.
\end{defn}

Special values are values of $]0,1[\times [-\overline{C}, \overline{C}]$ in which the optimal matching may abruptly change because of the presence of more than one pair of points with the same distance between abscissas (for $a \leq \frac{1}{2}$) or same distance between ordinates (for $a\geq \frac{1}{2}$). 
This discontinuity behaviour gives an obstruction to proving that the matching distance is realised only on vertical, horizontal and slope 1 lines. 
Indeed, the key for proving Theorem~\ref{minubyvtcrxeybt} is being able to continuously move in the space of parameters and not losing track of the points realising the optimal matching. 
When encountering a special value this continuity may be missing. 

\begin{figure}
    \centering
    \includegraphics[width=10cm]{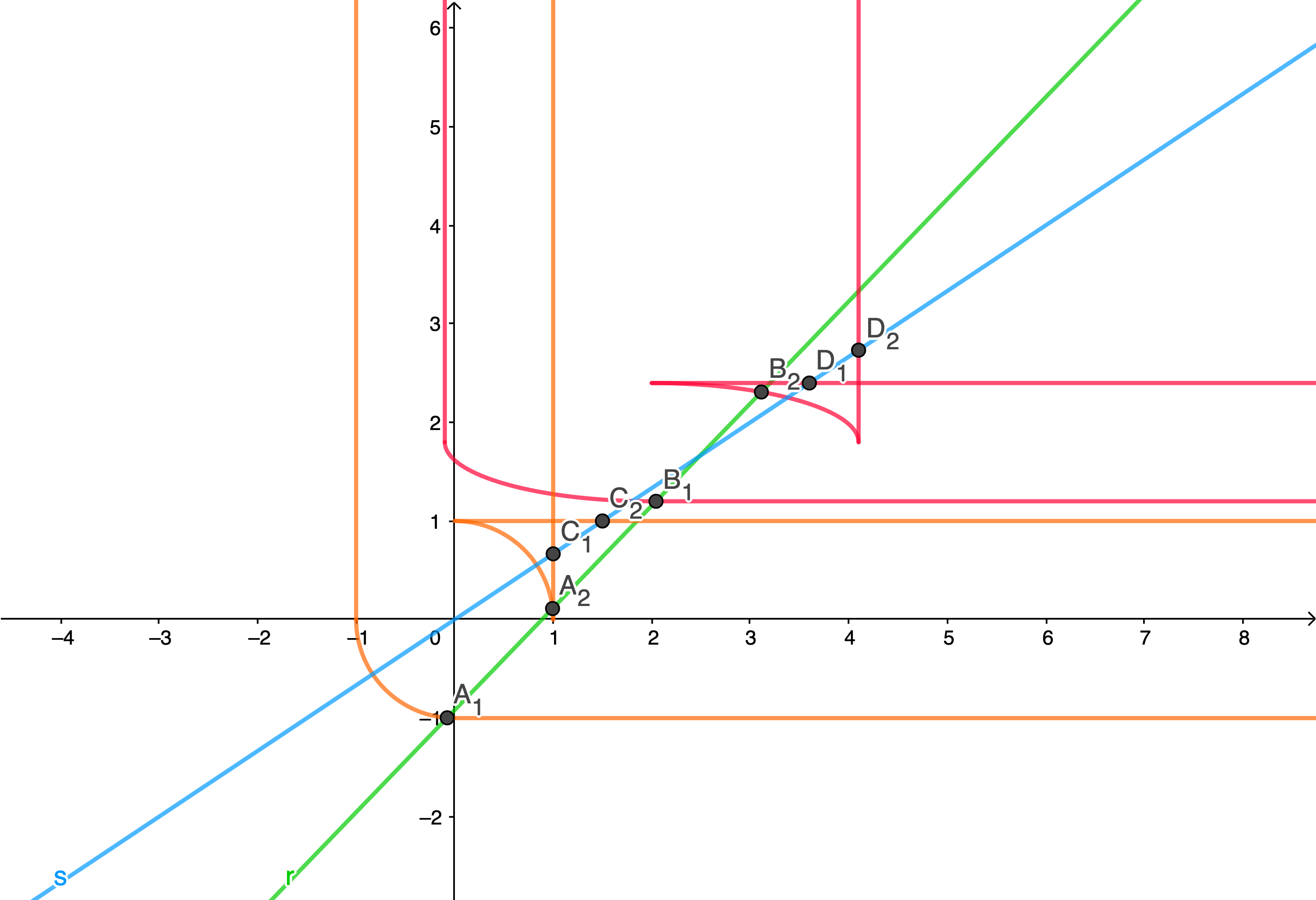}
    \caption{The light blue line corresponds to the pair $(0.6,0)$ and the green approximately to $(0.491,0.451)$. They are both special, since $\lvert x_{A_1}-x_{A_2}\rvert=\lvert x_{B_1}-x_{B_2}\rvert$ and $\lvert y_{C_1}-y_{C_2}\rvert=\lvert y_{D_1}-y_{D_2}\rvert$.}
    \label{special}
\end{figure}

Figure~\ref{special} shows two examples of lines associated with special values of $(f,g)$, with $f,g\colon S^2\to \mathbb{R}^2$, $f(x,y,z)=(x,z)$ and $g(x,y,z)=(2.1x+2, 0.6z+1.8)$. 
The green and light blue lines correspond respectively to the parameter values $(0.6,0)$ and $(0.491,0.451)$. 
The intersection points $A_1,A_2$ and $B_1, B_2$, between the green line and the extended Pareto grid have equal difference
between abscissas, thus $(0.6,0)$ is a special value.
On the other hand, the intersection points $C_1,C_2$ and $D_1, D_2$, between the light blue line
and the extended Pareto grid have equal difference between ordinates.
In particular, $(0.491,0.451)$ approximates a special value up to a $5\times 10^{-7}$ error.

\begin{prop}\label{vbhgjrueijlwks}
$\mathrm{Sp}(f,g)$ is closed in $]0,1[\times [-\overline{C}, \overline{C}]$. 
\end{prop}
\begin{proof}
First, we show that $\text{Sp}(f,g)\cap \left(]0,\frac{1}{2}]\times [-\overline{C}, \overline{C}]\right)$ is closed. 
Consider a sequence $\{(a_n, b_n)\}$ in $\text{Sp}(f,g)\cap \left(]0,\frac{1}{2}]\times [-\overline{C}, \overline{C}]\right)$ that converges to  $(\overline{a}, \overline{b})$ in $]0,\frac{1}{2}]\times [-\overline{C}, \overline{C}]$. Since such a sequence consists of special values of $(f,g)$, there exist two 
distinct sets $\{\alpha_p^n, \alpha_q^n\}$ and $\{\alpha_s^n, \alpha_t^n\}$ in $\text{Ctr}(f,g)$ such that
$c_1^n\lvert x_{P_n}-x_{Q_n}\rvert=c_2^n\lvert x_{S_n}-x_{T_n}\rvert$, where $P_n=r_{(a_n,b_n)}\cap \alpha_p^n$, $Q_n=r_{(a_n,b_n)}\cap \alpha_q^n$, 
$S_n=r_{(a_n,b_n)}\cap \alpha_s^n$ and $T_n=r_{(a_n,b_n)}\cap \alpha_t^n$ and $c_1^n, c_2^n\in\{1,2\}$, for every $n$. 
Since $\text{Ctr}(f,g)$ has finitely many contours, we can assume, up to subsequences, that the sequences $\{P_n\}$, $\{Q_n\}$, $\{S_n\}$ and $\{T_n\}$ lie respectively in the contours $\alpha_p^n=\alpha_p$, $\alpha_q^n=\alpha_q$, $\alpha_s^n=\alpha_s$ and $\alpha_t^n=\alpha_t$, for every $n$. 
For the same reason, we can assume that $c_1^n=c_1$ and $c_2^n=c_2$, for every $n$.
Since $\{(a_n,b_n)\}$ is convergent, it is also bounded. 
In particular, besides $-\overline{C}\le b_n\le \overline{C}$, 
there is $\overline{D}$ such that $0<\overline{D}\le a_n\le\frac{1}{2}$.
Then $\left(\bigcup_nr_{(a_n,b_n)}\right)\cap \left(\Gamma(f)\cup \Gamma(g)\right)$ is bounded below by the line $r_{(\frac{1}{2},\overline{C})}$ and above by the line $r_{(\overline{D},-\overline{C})}$.
Thus, $\{P_n\}$, $\{Q_n\}$, $\{S_n\}$ and $\{T_n\}$ converge, respectively, to $\overline{P}$, $\overline{Q}$, $\overline{S}$ and $\overline{T}$, up to restriction to subsequences. 
Since $c_1\lvert x_{P_n}-x_{Q_n}\rvert =c_2\lvert x_{S_n}-x_{T_n}\rvert$, their limits are also equal, so we have $c_1\lvert x_{\overline{P}}-x_{\overline{Q}}\rvert =c_2\lvert x_{\overline{S}}-x_{\overline{T}}\rvert$.
Since $\overline{P}$, $\overline{Q}$, $\overline{S}$ and $\overline{T}$ all lie in $r_{(\overline{a}, \overline{b})}$, $(\overline{a}, \overline{b})$ is also a special value of $(f,g)$, concluding that $\text{Sp}(f,g)\cap \left(]0,\frac{1}{2}]\times [-\overline{C}, \overline{C}]\right)$ is closed. 

Analogously, one can see that $\text{Sp}(f,g)\cap   \left([\frac{1}{2},1[\times [-\overline{C}, \overline{C}]\right)$ is closed. The set $\text{Sp}(f,g)$ is then a union of two closed sets, hence it is closed itself. 
\end{proof}

Let $\mathcal{S}$ be the set of all pairs $(\bar a,\bar b)$ in $[0,1] \times [-\overline{C},\overline{C}]$ realising the matching distance between $f$ and $g$, i.e., such that 
\[
D_{\mathrm{match}}(f,g)=d_B\left(\dgm\left(f_{(\bar a,\bar b)}^*\right),\dgm\left(g_{(\bar a,\bar b)}^*\right)\right).
\]
As observed about (\ref{contdb}), $d_B\left(\dgm\left(f_{(a,b)}^*\right),\dgm\left(g_{(a,b)}^*\right)\right)$ 
is a continuous function on 
$[0,1]\times[-\overline{C}, \overline{C}]$, thus it admits a 
maximum in its domain and $\mathcal{S}$ is not empty. 
Moreover, $\mathcal{S}$ is compact because it is the preimage of a point in $\mathbb{R}$ via a continuous function defined on a compact set.

Note that for any $(a,b)$ in $ [0,1] \times [-\overline{C}, \overline{C}]$, we have
\[
 d_B\left(\dgm\left(f_{( a, b)}^*\right),\dgm\left(g_{( a,b)}^*\right)\right) = \text{cost}(\sigma_{(a,b)}), \tag{5.1} \label{ciaociao}
\]
 where $\sigma_{(a,b)}$ is an optimal matching.
By applying a straightforward generalisation of Theorem 28 in~\cite{opt-matching} for arbitrary $n^{th}$
 persistence diagrams, one can see that such a matching always exists.
Theorem~\ref{ghjywuikdc} and the stability of the bottleneck distance with respect  to the uniform norm imply that $\text{cost}(\sigma_{(a,b)})$ can be seen as a continuous function in the variable $(a,b)$ in $[0,1] \times [-\overline{C}, \overline{C}]$.

\begin{defn}
Let $\sigma: \dgm_1\to\dgm_2$ be a matching between two persistence diagrams and let $X$ in $\dgm _1$ be such that $\mathrm{cost}(\sigma)=d(X,\sigma(X))$. The matching $\sigma$ is of type $(1)$ if $\Delta\notin \{X,\sigma(X)\}$, and of type $(2)$ if $\Delta\in \{X,\sigma(X)\}$.
\end{defn}

Observe that a matching can be both of type $(1)$ and type $(2)$.
We use this terminology in the proof of the following theorem.

\begin{thm}\label{minubyvtcrxeybt}
\[
\mathcal{S}\cap \left(\mathrm{Sp}(f,g) \cup \left(\left\{0,\frac{1}{2},1\right\}\times [-\overline{C},\overline{C}]\right) \right)\neq \emptyset.
\]
\end{thm}

\begin{figure}
    \centering
    \includegraphics[width=8cm]{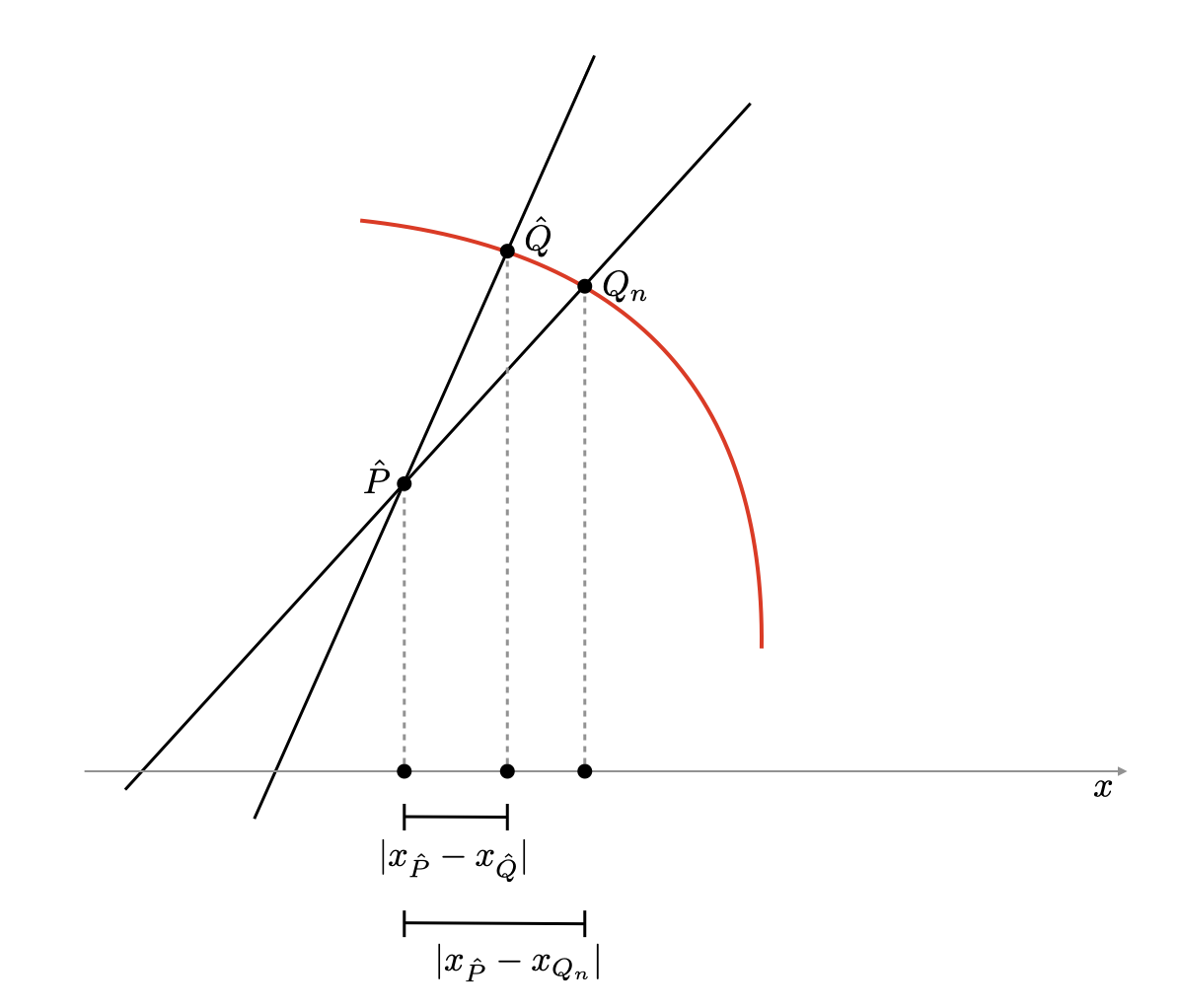}
    \caption{The clockwise rotation around $\widehat{P}$ increases the distance between the abscissas of the intersection points. This fact is used in the proof of Theorem~\ref{minubyvtcrxeybt} (case 1).}
    \label{proof_th}
\end{figure}

\begin{proof}
Assume by contradiction that every $(a,b)$ in $\mathcal{S}$ is not in $\text{Sp}(f,g)$ and that $a \neq 0,\frac{1}{2},1$.
Since $\mathcal{S}$ is compact, it is possible to take $(a,b)$ in $\mathcal{S}$ minimising the distance from the line $a=\frac{1}{2}$. 
Among these, consider $(\hat{a},\hat{b})$ and a corresponding matching $\widehat{\sigma} $ of minimum cost between $\dgm \left(f^*_{(\hat{a}, \hat{b})}\right)$ and $\dgm \left(g^*_{(\hat{a}, \hat{b})}\right)$.
If $0<\hat{a}< \frac{1}{2}$, the Position Theorem~\ref{posth} implies that there exist $\widehat\alpha$ and $\widehat\beta$ in $\text{Ctr}(f,g)$ 
intersecting $r_{(\hat{a},\hat{b})}$, such that $\widehat P=r_{(\hat a,\hat b)}\cap \widehat \alpha$ and $\widehat Q= r_{(\hat a,\hat b)}\cap\widehat \beta$ realise at least one of these properties: 
\begin{enumerate}
\item $\widehat P\in \Gamma(f)$, $\widehat Q\in \Gamma(g)$, and $D_{\mathrm{match}}(f,g)=\mathrm{cost} (\widehat\sigma)=\lvert x_{\widehat P}-x_{\widehat Q}\rvert$;
\item $\widehat P,\widehat Q\in \Gamma(f)$ or $\widehat P,\widehat Q\in \Gamma(g)$, and  $D_{\mathrm{match}}(f,g)=\mathrm{cost} (\widehat\sigma)=\frac{1}{2}\lvert x_{\widehat P}-x_{\widehat Q}\rvert$.
\end{enumerate}
Observe that the former matching is of type $(1)$ and the latter of type $(2)$. Note also that $x_{\widehat{P}}\neq x_{\widehat{Q}}$, and hence $\widehat{P}\neq \widehat{Q}$. If not, then $D_{\text{match}}(f,g)=0$, implying that any $(a,b)$ belongs to $\mathcal{S}$, including $\left(\frac{1}{2}, b\right)$, which is a contradiction.

Consider a sequence $\{(a_n, b_n)\}$ in  $]0,1[ \times [-\overline{C}, \overline{C}]$ such 
that these $(a_n, b_n)$ are chosen to identify lines obtained by rotating 
$r_{(a,b)}$ around $\widehat{P}$ clockwise in such a way that $(a_n, b_n) \to (\hat{a},\hat{b})$, where $\{a_n\}$ is a decreasing sequence.
Furthermore, given a sequence $\{\sigma_n\}$ of optimal matchings between $\dgm \left(f^*_{(a_n,b_n)}\right)$ and $\dgm \left(g
^*_{(a_n, b_n)}\right)$
we have that $\text{cost}(\sigma _n)\to\text{cost}(\widehat{\sigma})$  (see (\ref{ciaociao})).
Since $\text{Sp}(f,g) \cup \left(\left\{0,\frac{1}{2},1\right\}\times [-\overline{C},\overline{C}]\right)$ is closed, by Proposition~\ref{vbhgjrueijlwks}, and $(\hat a,\hat b)$ does not belong to this set, we can assume that the sequence $\{(a_n,b_n)\}$ also has no points in this set.
Hence, for any $n$ in $\mathbb{N}$ there exists a pair $\{P_n, Q_n\}$ in $r_{(a_n, b_n)}\cap(\Gamma(f)\cup \Gamma(g))$ for which at least one of the following properties holds:
\begin{enumerate}[(A)]
\item $P_n\in \Gamma(f)$, $ Q_n\in \Gamma(g)$ and $\mathrm{cost} (\sigma_n)=\lvert x_{ P_n}-x_{ Q_n}\rvert$;
\item $P_n,Q_n\in \Gamma(f)$ or $P_n,Q_n\in \Gamma(g)$, and  $\mathrm{cost} (\sigma_n)=\frac{1}{2}\lvert x_{P_n}-x_{Q_n}\rvert$.
\end{enumerate}

Up to subsequences, we can assume that the matchings $\sigma_n$ are either all of type $(1)$ or all of type $(2)$.
We now show that $\{P_n\}$ and $\widehat{P}$ belong to the same contour in $\text{Ctr}(f,g)$, and $\{Q_n\}$ and $\widehat{Q}$ also belong to the same contour in $\text{Ctr}(f,g)$.  
Analogously to the proof of Proposition~\ref{vbhgjrueijlwks} we may observe that
the set $\left(\bigcup_n r_{(a_n,b_n)}\right) \cap(\Gamma(f)\cup \Gamma(g))$
is a bounded subset of $\Gamma(f)\cup \Gamma(g)$.
Thus, $\{P_n\}$ and $ \{Q_n\}$ are convergent up to subsequences in the closed set $\Gamma(f)\cup \Gamma(g)$, respectively, to $\overline{P}$ and $\overline{Q}$.
By assumption, there are only a finite number of contours, thus there exists at least a contour in $\text{Ctr}(f,g)$ for each sequence, $\{P_n\}$ and $\{Q_n\}$, containing infinitely many points of the sequence.
Hence, we can assume that each sequence, up to subsequences, lies entirely on a single contour in $\text{Ctr}(f,g)$, i.e., 
we can suppose that for every $n$ in $\mathbb{N}$, $P_n$ is in $\overline{\alpha}$ and $Q_n$ is in $\overline{\beta}$, with $\overline{\alpha}$ and $\overline{\beta}$ in $\text{Ctr}(f,g)$.
Since contours are closed, $\overline{P}$ belongs to $\overline{\alpha}$ and $\overline{Q}$ belongs to $\overline{\beta}$. We observe that $\{\overline{P}, \overline{Q}\} \subseteq r_{(\hat{a},\hat{b})}$. Furthermore,  we have that
\[c'\lvert x_{\widehat{P}}-x_{\widehat{Q}}\rvert = \text{cost}(\widehat{\sigma}) = \lim_{n \to \infty} \text{cost}(\sigma _n)= \lim_{ n \to \infty} c''\lvert x_{P_n}-x_{Q_n}\rvert=c''\lvert x_{\overline{P}}-x_{\overline{Q}}\rvert
\]
where $c',c''$ in $\{\frac{1}{2},1\}$.
If $\{\widehat{\alpha}, \widehat{\beta}\}\neq \{\overline{\alpha}, \overline{\beta}\}$, then $(\hat{a}, \hat{b})$ is a special value, contradicting the initial assumption. 
Thus, $\{\widehat{\alpha}, \widehat{\beta}\}= \{\overline{\alpha}, \overline{\beta}\}$. 
Without loss of generality, by possibly exchanging the roles of the contours $\widehat \alpha$ and $\widehat\beta$, and of
the points $\widehat P$ and $\widehat Q$, we can assume that $\widehat{\alpha}=\overline{\alpha}$,
$\widehat{\beta}=\overline{\beta}$, $\widehat{P}=\overline{P}$ and $\widehat{Q}=\overline{Q}$.
Consequently, by the fact that $\{P_n\}$ and $\widehat{P}$ are contained in the same line $r_{(a_n,b_n)}$ and the same contour $\widehat\alpha$, $P_n=\widehat{P}$ for every $n$, since a contour and a positive slope line can meet in at most one point.  

\textbf{Case 1.} Assume that $\widehat \sigma$ and $\sigma_n$ are both of the same type for every $n$.
Since $Q_n$ belongs to $\overline{\beta}$ in $\text{Ctr}(f,g)$ for any $n$, one can easily check that $\lvert x_{\widehat{P}}-x_{Q_n}\rvert\ge\lvert x_{\widehat{P}}-x_{\widehat{Q}}\rvert$ (see Figure~\ref{proof_th}),
and hence $\text{cost}(\sigma_n)\ge \text{cost}(\widehat{\sigma})$. If the equality holds there is a  contradiction with the assumption of $(\hat{a}, \hat{b})$ minimising the distance from the line $a=\frac{1}{2}$, since $\left\lvert a_n- \frac{1}{2}\right\rvert <\left\lvert \hat a- \frac{1}{2}\right\rvert$. 
If the strict inequality holds, there is a contradiction with the assumption of $\widehat{\sigma}$ being in $\mathcal{S}$. 

\textbf{Case 2.} Assume that all $\sigma_n$ and $\widehat \sigma$ are of different types. This means that $\text{cost}(\sigma _n)=c'\lvert x_{\widehat{P}}-x_{Q_n}\rvert$, $\text{cost}(\widehat{\sigma})=c''\lvert x_{\widehat{P}}-x_{\widehat{Q}}\rvert$, with $c'\neq c''$ and $c',c''$ in $\{\frac{1}{2},1\}$, and $ c' \lvert x_{\widehat{P}}-x_{Q_n}\rvert \to c''\lvert x_{\widehat{P}}-x_{\widehat{Q}}\rvert$. However, since $Q_n\to \widehat{Q}$, $c'\lvert x_{\widehat{P}}-x_{Q_n}\rvert\to c'\lvert x_{\widehat{P}}-x_{\widehat{Q}}\rvert$. Thus, $c''\lvert x_{\widehat{P}}-x_{\widehat{Q}}\rvert=c'\lvert x_{\widehat{P}}-x_{\widehat{Q}}\rvert$, which is a contradiction since $D_{\text{match}}(f,g)\neq 0$ and, hence, $x_{\widehat{P}}\neq x_{\widehat{Q}}$.

Inverting the role of abscissas and ordinates as described by the Position Theorem~\ref{posth} and rotating the lines counterclockwise, one can see that an analogous procedure holds for $\frac{1}{2}<\hat a<1$. 
\end{proof}

\section{Conclusions}
In this article we took advantage of the differential structure associated with smooth functions from a Riemannian manifold $M$ to $\mathbb{R}^2$ to characterise some geometric properties of the matching distance. 
We proved that the filtering lines that actually contribute to the computation of the matching distance are horizontal, vertical, of slope 1, or they are associated with parameter values in the special set. 
This new approach to the computation of the matching distance could lead to new effective algorithms. 
In this direction, we would like to highlight an open question that arose during our work. 
We have not yet provided a characterisation of the special set. 
However, we conjecture that the special set consists of a collection of curves, up to a small perturbation of the filtering functions. 

\begin{figure}
\centering
\includegraphics[width=7cm]{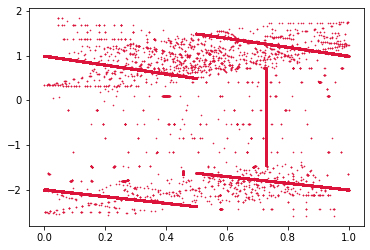}
\caption{Approximation of a special set.}
\label{newest_tol}
\end{figure}

Figure~\ref{newest_tol} shows a selection of points in the special set for the functions $f,g\colon S^2\to\mathbb{R}^2$, where $S^2=\{(x,y,z)\mid x^2+y^2+z^2=1\}$, $f(x,y,z)=(x+1, z-1)$ and $g(x,y,z)=(0.75x-2,0.75z+2)$.
One may notice clear segments, two of which, on the left, correspond to values identifying lines through intersections of contours. 
Such lines are in fact always associated with special values. 

\bibliographystyle{plain}
\bibliography{references}

\begin{thebibliography}{10}

\bibitem{comp_disc}
Asilata Bapat, Robyn Brooks, Celia Hacker, Claudia Landi, Barbara~I. Mahler,
  and Elizabeth~R. Stephenson.
\newblock Computing the matching distance of 2-parameter persistence.
\newblock {\em arXiv:2210.12868}.

\bibitem{comp2}
Silvia Biasotti, Andrea Cerri, Patrizio Frosini, and Daniela Giorgi.
\newblock A new algorithm for computing the 2-dimensional matching distance
  between size functions.
\newblock {\em Pattern Recognition Letters}, 32(14):1735--1746, 2011.

\bibitem{comp_asymp}
Havard Bjerkevik and Michael Kerber.
\newblock Asymptotic improvements on the exact matching distance for
  2-parameter persistence.
\newblock {\em arXiv: 2111.10303}.

\bibitem{pers_Betti_stable}
Andrea Cerri, Barbara Di~Fabio, Massimo Ferri, Patrizio Frosini, and Claudia
  Landi.
\newblock Betti numbers in multidimensional persistent homology are stable
  functions.
\newblock {\em Math. Methods Appl. Sci.}, 36(12):1543--1557, 2013.

\bibitem{coher_match}
Andrea Cerri, Marc Ethier, and Patrizio Frosini.
\newblock On the geometrical properties of the coherent matching distance in
  2{D} persistent homology.
\newblock {\em J. Appl. Comput. Topol.}, 3(4):381--422, 2019.

\bibitem{comp1}
Andrea Cerri and Patrizio Frosini.
\newblock A new approximation algorithm for the matching distance in
  multidimensional persistence.
\newblock {\em Journal of Computational Mathematics}, 38(2):291--309, 2020.

\bibitem{filtered-persistence}
Herbert Edelsbrunner and Dmitriy Morozov.
\newblock Persistent homology: theory and practice.
\newblock In {\em European {C}ongress of {M}athematics}, pages 31--50. Eur.
  Math. Soc., Z\"{u}rich, 2013.

\bibitem{comp_exact}
Michael Kerber, Michael Lesnick, and Steve Oudot.
\newblock Exact computation of the matching distance on 2-parameter persistence
  modules.
\newblock In {\em 35th {I}nternational {S}ymposium on {C}omputational
  {G}eometry}, volume 129 of {\em LIPIcs. Leibniz Int. Proc. Inform.}, pages
  Art. No. 46, 15. Schloss Dagstuhl. Leibniz-Zent. Inform., Wadern, 2019.

\bibitem{efficient-approx}
Michael Kerber and Arnur Nigmetov.
\newblock Efficient approximation of the matching distance for 2-parameter
  persistence.
\newblock In {\em 36th {I}nternational {S}ymposium on {C}omputational
  {G}eometry}, volume 164 of {\em LIPIcs. Leibniz Int. Proc. Inform.}, pages
  Art. No. 53, 16. Schloss Dagstuhl. Leibniz-Zent. Inform., Wadern, 2020.

\bibitem{inter}
Michael Lesnick.
\newblock The theory of the interleaving distance on multidimensional
  persistence modules.
\newblock {\em Found. Comput. Math.}, 15(3):613--650, 2015.

\bibitem{wan}
Y.~H. Wan.
\newblock Morse theory for two functions.
\newblock {\em Topology}, 14(3):217--228, 1975.

\end{thebibliography}

\end{document}